\documentclass{amsart}

\usepackage{amsmath,amsthm,amsfonts,amssymb,amscd,latexsym}

\newtheorem{thm}{Theorem}[section]
\newtheorem{lem}[thm]{Lemma}
\newtheorem{pro}[thm]{Proposition}
\newtheorem{cor}[thm]{Corollary}

\numberwithin{equation}{section}

\newcommand{\id}{\mathrm{id}}
\newcommand{\Hom}{\mathrm{Hom}}
\newcommand{\End}{\mathrm{End}}
\newcommand{\ann}{\mathrm{Ann}}
\newcommand{\triv}{\mathrm{triv}}
\newcommand{\rad}{\mathrm{Rad}}
\newcommand{\leib}{\mathrm{Leib}}
\newcommand{\lie}{\mathrm{Lie}}
\newcommand{\der}{\mathrm{Der}}

\newcommand{\HL}{\mathrm{HL}}
\newcommand{\ad}{\mathrm{ad}}
\newcommand{\rk}{\mathrm{rk}}

\newcommand{\Z}{\mathbb{Z}}
\newcommand{\C}{\mathbb{C}}
\newcommand{\F}{\mathbb{F}}

\newcommand{\lf}{\mathfrak{L}}
\newcommand{\nf}{\mathfrak{N}}
\newcommand{\af}{\mathfrak{A}}
\newcommand{\I}{\mathfrak{I}}
\newcommand{\g}{\mathfrak{g}}
\newcommand{\s}{\mathfrak{s}}

\newcommand{\slf}{\mathfrak{sl}}

\begin{document}


\title[Semi-simple Leibniz algebras]{Semi-simple Leibniz algebras I}

\author{J\"org Feldvoss}
\address{Department of Mathematics and Statistics, University of South
Alabama, Mobile, AL 36688-0002, USA}
\email{jfeldvoss@southalabama.edu}



\begin{abstract}
The goal of this paper is to describe the structure of finite-dimen\-sional
semi-simple Leibniz algebras in characteristic zero. Our main tool in this
endeavor are hemi-semidirect products. One of the major results of this
paper is a simplicity criterion for hemi-semidirect products. In addition, we
characterize when a hemi-semidirect product is semi-simple or Lie-simple.
Using these results we reduce the classification of finite-dimensional
semi-simple Leibniz algebras over fields of characteristic zero to the
well-known classification of finite-dimensional semi-simple Lie algebras
and their finite-dimensional irreducible modules. As one consequence of our
structure theorem, we determine the derivation algebra of a finite-dimensional
semi-simple Leibniz algebra in characteristic zero as a vector space. This
generalizes a recent result of Ayupov et al.\ from the complex numbers
to arbitrary fields of characteristic zero.
\end{abstract}


\date{January 10, 2024}
          
\maketitle


\section*{Introduction}


In this paper we investigate the structure of semi-simple
Leibniz algebras. We mostly consider finite-dimensional
algebras over fields of characteristic zero, although the
results of Section \ref{hemisemi} are valid more generally
for arbitrary ground fields and algebras of arbitrary dimension.

In the first section we recall some basic definitions for Leibniz
algebras and Leibniz (bi)modules. In particular, we define
(left) hemi-semidirect products which will be crucial for our
main results in the other sections. Hemi-semidirect products
have been introduced by Kinyon and Weinstein in their study
of reductive homogeneous spaces (see \cite{KW}), and it
turns out that many Leibniz algebras are hemi-semidirect
products. For later use we determine the Leibniz kernel and
the derived subalgebra as well as the left and right center
of an arbitrary hemi-semidirect product.

The main results of the paper are contained in Sections
\ref{hemisemi} and \ref{structhm}. In Section~\ref{hemisemi}
we characterize the semi-simplicity or (Lie-)simplicity of a
hemi-semidirect product (see Theorems \ref{semisim},
\ref{liesim}, and \ref{sim}). In part, the missing hypothesis
in Example 2.3 of \cite{ORT} and in Example 1 at the end of
Section 2 of \cite{GVKO} as well as the lack of any proof have
been the motivation for writing this paper, and especially, for
proving the simplicity criterion Theorem \ref{sim}. Moreover, as
far as I can see, most of the existing papers on (semi-)simple
Leibniz algebras (for example, \cite{ORT}, \cite{GVKO},
\cite{CGVOK}, \cite{COK}, \cite{AKOZ}) do not use a conceptual
approach, and the latter is exactly what we will present here.
In particular, in Section \ref{hemisemi} we discuss how our
results can be employed to prove some of the statements in previous
papers on (semi-)simple Leibniz algebras or sometimes how to
complete the proofs of these statements.

The results of Section \ref{hemisemi} enable us then in Section
\ref{structhm} to obtain a structural characterization of
finite-dimensional semi-simple and (Lie-)simple (left) Leibniz
algebras in characteristic zero as certain (left) hemi-semidirect
products (see Theorems~\ref{semisimchar0}, \ref{liesimchar0},
and \ref{simchar0}). Since for fields of characteristic
zero the classification of finite-dimensional semi-simple Lie
algebras and their finite-dimensional irreducible modules is
well-known, this allows us to construct all finite-dimensional
semi-simple and (Lie-)simple Leibniz algebras over fields of
characteristic zero. Note that the proofs in this section
crucially use Levi's theorem for Leibniz algebras (see
\cite[Proposition~2.4]{P} or \cite[Theorem 1]{B}) and
Weyl's theorem on complete reducibility (see \cite[Theorem
6.3]{H}), and therefore one cannot expect similar results
for infinite-dimensional Leibniz algebras or finite-dimensional
Leibniz algebras over fields of non-zero characterstics. But
even in these cases one can say quite a bit, and we will
consider these matters in the forthcoming second part of
this paper. We conclude this section by briefly discussing
the fact that, contrary to Lie algebras, finite-dimensional
semi-simple Leibniz algebras in characteristic zero are not
necessarily a direct sum of (Lie-)simple Leibniz algebras
(see also \cite[Examples 2 and 3]{GVKO}) and how our
structure theorems can be used to construct examples
showing this fact in a conceptual manner.

In the last section we give some applications of our structure
theorem for finite-dimensional semi-simple Leibniz algebras. In
particular, we show that the right center of such an algebra
vanishes, and as a consequence, we obtain that every
finite-dimensional semi-simple left central Leibniz algebra in
characteristic zero is already a Lie algebra (see Proposition
\ref{rightcent} and Corollary \ref{symmsemisim}). Moreover,
we can describe the derivation algebra of a finite-dimensional
semi-simple Leibniz algebra as a vector space. This generalizes
one of the main results in a paper of Ayupov, Kudaybergenov,
Omirov, and Zhao on the dimension of such a derivation algebra
(see \cite[Theorem 4.5\,(c)]{AKOZ}) from the complex numbers to
an arbitrary field of characteristic zero (see Theorem \ref{semisimder}).

In this paper we will follow the notation used in \cite{F} and
\cite{FW}. An algebra without any specification will be a vector
space with a bilinear multiplication that not necessarily satisfies
any other identities. Ideals will always be two-sided ideals. All
vector spaces, (bi)modules, and algebras are defined over an
arbitrary field which is only explicitly mentioned when some
additional assumptions on the ground field are made or this
enhances the understanding of the reader. For a subset
$X$ of a vector space $V$ over a field $\F$ we let $\langle
X\rangle_\F$ be the subspace of $V$ spanned by $X$.
Moreover, we will denote the space of linear transformations
from an $\F$-vector space $V$ to an $\F$-vector space $W$
by $\Hom_\F(V,W)$. In particular, $\End_\F(V):=\Hom_\F(V,V)$
is the space of linear operators on $V$. Finally, the identity
function on a set $X$ will be denoted by $\id_X$, the set of
integers will be denoted by $\Z$, and $\Z_2=\Z/2\Z=
\{\overline{0},\overline{1}\}$ is the additive group of integers
modulo 2.


\section{Preliminaries}\label{prelim}


In this section we collect some definitions and basic results that will be useful in the
remainder of the paper, and for which the reader is referred to
this section without explicitly mentioning all the time the exact
place. Moreover, we determine the Leibniz kernel and the derived
subalgebra as well as the left and right center of a left hemi-semidirect
product.

A {\em left Leibniz algebra\/} is an algebra $\lf$ such that every left
multiplication operator $L_x:\lf\to\lf$, $y\mapsto xy$ is a derivation.
This is equivalent to the identity
$$x(yz)=(xy)z+y(xz)$$
for all $x,y,z\in\lf$, which in turn is equivalent to the identity
$$(xy)z=x(yz)-y(xz)$$
for all $x,y,z\in\lf$. There is a similar definition of a {\em right Leibniz
algebra\/}, but in this paper we will mostly consider left Leibniz algebras.
Following Mason and Yamskulna \cite[pp.\ 1/2]{MY} we say that a left
and right Leibniz algebra is {\em symmetric\/}.

Note that every Lie algebra is a symmetric Leibniz algebra. On the
other hand, every Leibniz algebra has an important ideal, its Leibniz
kernel, that measures how much the Leibniz algebra deviates from
being a Lie algebra. Namely, let $\lf$ be a Leibniz algebra over a
field $\F$. Then
$$\leib(\lf):=\langle x^2\mid x\in\lf\rangle_\mathbb{F}$$
is called the {\em Leibniz kernel\/} of $\lf$. Note that the Leibniz
kernel of a Leibniz algebra $\lf$ is contained in the {\em derived
subalgebra\/}
$$\lf^2:=\langle xy\mid x,y\in\lf\rangle_\F$$
of $\lf$. The Leibniz kernel $\leib(\lf)$ is an abelian ideal of $\lf$,
and $\leib(\lf)\ne\lf$ whenever $\lf\ne 0$ (see \cite[Proposition
2.20]{F}). Moreover, $\lf$ is a Lie algebra exactly when $\leib
(\lf)=0$. By definition of the Leibniz kernel, $\lf_\lie:=\lf/\leib(\lf)$
is a Lie algebra which we call the {\em canonical Lie algebra\/}
associated to $\lf$. Let
$$C^\ell(\lf):=\{c\in\lf\mid\forall\,x\in\lf:cx=0\}$$
denote the {\em left center\/} of a left Leibniz algebra $\lf$. Then
$C^\ell(\lf)$ is an abelian ideal of $\lf$ (see \cite[Proposition 2.8]{F})
which contains the Leibniz kernel $\leib(\lf)$ (see \cite[Proposition
2.13]{F}). In Section \ref{appl} we will also need the {\em right
center\/}
$$C^r(\lf):=\{c\in\lf\mid\forall\,x\in\lf:xc=0\}$$
of a left Leibniz algebra $\lf$.

Next, we briefly discuss left modules and bimodules over left Leibniz
algebras. Let $\lf$ be a left Leibniz algebra over a field $\F$. A {\em
left $\lf$-module\/} is a vector space $M$ over $\F$ with an $\F$-bilinear
left $\lf$-action $\lf\times M\to M$, $(x,m)\mapsto x\cdot m$ such that
$$(xy)\cdot m=x\cdot(y\cdot m)-y\cdot(x\cdot m)$$
is satisfied for every $m\in M$ and all $x,y\in\lf$.

Note that the correct concept of a module for a left Leibniz algebra $\lf$
is the notion of a Leibniz bimodule (see \cite[Section 3]{F} for a motivation
behind this definition, namely, that the {\em semidirect product\/} $\lf
\ltimes M$ of a left Leibniz algebra $\lf$ and a Leibniz $\lf$-bimodule $M$
with underlying vector space $\lf\times M$ and multiplication
$$(g,m)(g^\prime,m^\prime):=(gg^\prime,g\cdot m^\prime+m\cdot g^\prime)$$
for any elements $g,g^\prime\in\g$ and $m,m^\prime\in M$, is a left Leibniz
algebra). An {\em $\lf$-bimodule\/} is a vector space $M$ with an $\F$-bilinear
left $\lf$-action and an $\F$-bilinear right $\lf$-action that satisfy the following
compatibility conditions:
\begin{enumerate}
\item[(LLM)] \hspace{2.5cm}$(xy)\cdot m=x\cdot(y\cdot m)-y\cdot(x\cdot m)$
\item[(LML)] \hspace{2.5cm}$(x\cdot m)\cdot y=x\cdot(m\cdot y)-m\cdot(xy)$
\item[(MLL)] \hspace{2.5cm}$(m\cdot x)\cdot y=m\cdot(xy)-x\cdot(m\cdot y)$
\end{enumerate}
for every $m\in M$ and all $x,y\in\lf$. It is an immediate consequence of
(LLM) that every Leibniz bimodule is a left Leibniz module.

Since (LLM) implies that $\leib(\lf)M=0$, every left $\lf$-module is a left
module over the canonical Lie algebra $\lf_\lie=\lf/\leib(\lf)$ of $\lf$ in a
natural way, and vice versa. Consequently, many properties of left Leibniz
modules follow from the corresponding properties of left modules for the
canonical Lie algebra.

The usual definitions of the notions of {\em sub{\rm (}\hspace{-.5mm}bi{\rm
)}\hspace{-.5mm}module\/}, {\em irreducibility\/}, {\em complete reducibility\/},
etc., hold for left Leibniz modules and Leibniz bimodules. Note that by
definition an irreducible (bi)module is always non-zero, but a completely
reducible (bi)module is allowed to be zero.

Finally, an $\lf$-bimodule $M$ is called {\em trivial\/} if both its left and
right $\lf$-actions are trivial (i.e., $x\cdot m=0=m\cdot x$ for all $x\in
\lf$ and $m\in M$), and when only its right $\lf$-action is trivial, then
$M$ is said to be {\em anti-symmetric\/}.
\vspace{.2cm}

{\bf Example:} Every left Leibniz algebra $\lf$ is a bimodule over
itself via its left and right multiplication, the so-called {\em adjoint
bimodule\/} $\lf_\ad$ of $\lf$, and as $\leib(\lf)\lf=0$, it follows
that $\leib(\lf)$ is an anti-symmetric $\lf$-subbimodule of $\lf_\ad$.
\vspace{.2cm}

A very important tool to construct Leibniz algebras from Lie algebras
and Lie modules are hemi-semidirect products, and indeed, the latter
will be crucial in this paper. As already mentioned in the introduction,
hemi-semidirect products have been defined for the first time by
Kinyon and Weinstein (see \cite[Example 2.2]{KW}), but are also,
at least implicitly, employed elsewhere (see \cite[p.\ 1049]{DA},
\cite[Example~2.3]{ORT}, \cite[Example 1]{GVKO}, \cite[Section 3,
(7)]{MY}, \cite[Section 2, (1)]{COK}, \cite[Definition 2.8]{AKOZ},
and \cite[Definition 1.5]{OW}).

The {\em left hemi-semidirect product\/} $\g\ltimes_\ell M$ of a
Lie algebra $\g$ and a left $\g$-module $M$ is defined as follows.
(Note that a left module over a Lie algebra satisfies the same defining
identity as a left module over a left Leibniz algebra, see (LLM).) 
The underlying vector space of $\g\ltimes_\ell M$ is the Cartesian
product $\g\times M$ of $\g$ and $M$ with componentwise addition
and componentwise scalar multiplication, and the multiplication on
$\g\ltimes_\ell M$ is given by
$$(g,m)(g^\prime,m^\prime):=([g,g^\prime],g\cdot m^\prime)$$
for any elements $g,g^\prime\in\g$ and $m,m^\prime\in M$, where
$[-,-]$ denotes the Lie bracket on $\g$ and $\cdot$ denotes the left
action of $\g$ on $M$.\footnote{As usual, for any subset $S$
of $\g$ we will identify $S\times\{0\}$ with $S$, and similarly, for
any subset $T$ of $M$ we will identify $\{0\}\times T$ with $T$.}

It is known that $\g\ltimes_\ell M$ is a left Leibniz algebra (see
\cite[Example 2.2]{KW}). This can also be verified by the computations
on p.\ 128 in \cite{F} as $M$ can be made into an anti-symmetric
$\g$-bimodule $M_a$, the so-called {\em anti-symmetrization\/}
of $M$ (see \cite[Proposition 3.15\,(a)]{F}). Namely, then the left
hemi-semidirect product $\g\ltimes_\ell M$ is just the semidirect
product $\g\ltimes M_a$ of $\g$ and the anti-symmetrization
$M_a$ of $M$.
The smallest non-trivial example is the left hemi-semidirect product
$\af=\F h\ltimes_\ell\F e$ of a one-dimensional Lie algebra $\F h$
and a one-dimensional module $\F e$ with action $h\cdot e=e$
(see \cite[Example~2.3]{F} for an explicit definition of the multiplication
of $\af$).

In Sections \ref{hemisemi} and \ref{appl} we will need the
following descriptions of the Leibniz kernel as well as the left
and right center of a left hemi-semidirect product, respectively.
In addition, we determine the derived subalgebra of a left
hemi-semidirect product which is closely related to its Leibniz kernel.
The first part was already observed by Mason and Yamskulna (see
\cite[Section 3, p.\ 3]{MY}). In order to be able to formulate Lemma
\ref{prophemisemi}, we need some more notation, namely,
the {\em center\/}
$$C(\g):=C^\ell(\g)=C^r(\g)$$
of a Lie algebra $\g$ as well as the {\em $\g$-annihilator\/}
$$\ann_\g(M):=\{g\in\g\mid\forall\,m\in M:g\cdot m=0\}\,,$$
the {\em space of $\g$-invariants\/}
$$M^\g:=\{m\in M\mid\forall\,g\in\g:g\cdot m=0\}\,,$$
and the subspace
$$\g M:=\langle g\cdot m\mid g\in\g,m\in M\rangle_\F$$
of a left $\g$-module $M$.

\begin{lem}\label{prophemisemi}
Let $\g$ be a Lie algebra, and let $M$ be a left $\g$-module. Then
the following statements hold:
\begin{enumerate}
\item[(a)] $\leib(\g\ltimes_\ell M)=\g M$,
\item[(b)] $(\g\ltimes_\ell M)^2=\g^2\times \g M$,
\item[(c)] $C^\ell(\g\ltimes_\ell M)=[C(\g)\cap\ann_\g(M)]\times M$,
\item[(d)] $C^r(\g\ltimes_\ell M)=C(\g)\times M^\g$.
\end{enumerate}
\end{lem}

\begin{proof}
Let $g\in\g$ and $m\in M$ be arbitrary elements. In order to prove
(a), we consider the identity
$$(g,m)(g,m)=([g,g],g\cdot m)=(0,g\cdot m)\,.$$
Since $\leib(\g\ltimes_\ell M)$ and $\g M$ are subspaces of $\g
\ltimes_\ell M$ and $M$, respectively, this shows the assertion.

For the proofs of the remaining statements consider the identity
$$(a,m)(b,n)=([a,b],a\cdot n)$$
for any elements $a,b\in\g$ and $m,n\in M$.
\end{proof}

Lemma \ref{prophemisemi} can be illustrated by the hemi-semidirect
product $\af=\F h\ltimes_\ell\F e$ mentioned earlier (see \cite[Examples
2.10 and 2.16]{F}).


\section{(Semi-)Simplicity of Hemi-semidirect products}\label{hemisemi}


A Leibniz algebra $\lf$ is called {\em semi-simple\/} if $\leib(\lf)$
contains every solvable ideal of $\lf$ (see \cite[Section 7]{F}).
Note that for this definition we do not need to assume that $\lf$
satisfies any finiteness conditions. In particular, it follows that
a finite-dimensional Leibniz algebra $\lf$ is semi-simple exactly
when $\leib(\lf)=\rad(\lf)$, where $\rad(\lf)$ denotes the
{\em radical\/} (= largest solvable ideal) of $\lf$ (see
\cite[Proposition 7.4]{F}). Usually this is employed as the
definition for semi-simplicity (see \cite[Definition 2.6]{GVKO}
and \cite[Definition 5.2]{DMS}). Note also that the
semi-simplicity of a Leibniz algebra $\lf$ is equivalent to
the semi-simplicity of its canonical Lie algebra $\lf_\lie=
\lf/\leib(\lf)$ (see \cite[Proposition~7.8]{F}). Therefore,
in the literature sometimes the semi-simplicity of the canonical
Lie algebra is used as the definition for the semi-simplicity
of a Leibniz algebra (for example, see \cite[Definition 2.9]{AKOZ}).

Our first result characterizes when a left hemi-semidirect product
is semi-simple.

\begin{thm}\label{semisim}
Let $\g$ be a Lie algebra, and let $M$ be a left $\g$-module. 
Then the left hemi-semidirect product $\g\ltimes_\ell M$ is
semi-simple if, and only if, $\g$ is semi-simple and $\g M=M$.
\end{thm}

\begin{proof}
Let $\lf:=\g\ltimes_\ell M$ denote the left hemi-semidirect product
of $\g$ and $M$.  Suppose first that $\lf$ is semi-simple. It is clear
from the definition of the hemi-semidirect product that $M$ is an
abelian ideal of $\lf$, and therefore we obtain that $M\subseteq
\leib(\lf)$. On the other hand, we deduce from Lemma
\ref{prophemisemi}\,(a) that $\leib(\lf)=\g M\subseteq M$, and
thus $\leib(\lf)=M$. Moreover, it follows from the semi-simplicity
of $\lf$ and \cite[Proposition 7.8]{F} that
$$\g\cong\lf/M=\lf/\leib(\lf)=\lf_\lie$$
is semi-simple.

Conversely, suppose that $\g$ is semi-simple and $\g M=M$. By
virtue of Lemma~\ref{prophemisemi} \,(a), we have that $\leib(\lf)
=\g M$, and therefore $\leib(\lf)=M$. Since
$$\lf_\lie=\lf/\leib(\lf)=\lf/M\cong\g\,,$$
and by hypothesis $\g$ is semi-simple, the assertion follows from
Proposition 7.8 in \cite{F}.
\end{proof}

Let us mention the following application of Theorem \ref{semisim}.
Namely, it is a consequence of (the analog of) Theorem \ref{semisim}
(for right Leibniz algebras) that the Leibniz algebras in \cite[Section
4]{CGVOK} are semi-simple.

A Leibniz algebra $\lf$ is called {\em Lie-simple\/} if its canonical
Lie algebra $\lf_\lie=\lf/\leib(\lf)$ is simple (see \cite[Section 1,
p.\ 527]{GVKO}). The next result is the direct analogue of Theorem
\ref{semisim} and characterizes when a left hemi-semidirect product
is Lie-simple.

\begin{thm}\label{liesim}
Let $\g$ be a Lie algebra, and let $M$ be left $\g$-module. Then
the left hemi-semidirect product $\g\ltimes_\ell M$ is Lie-simple if,
and only if, $\g$ is simple and $\g M=M$.
\end{thm}

\begin{proof}
Let $\lf:=\g\ltimes_\ell M$ denote the left hemi-semidirect product
of $\g$ and $M$. Suppose first that $\lf:=\g\ltimes_\ell M$ is Lie-simple.
Since Lie-simple Leibniz algebras are semi-simple (see \cite[Proposition
7.8]{F}), we can argue as in the proof of Theorem~\ref{semisim}.

Conversely, suppose that $\g$ is simple and $\g M=M$. As in the
proof of Theorem~\ref{semisim}, we have that $\leib(\lf)=\g M=M$.
Since by hypothesis
$$\lf_\lie=\lf/\leib(\lf)=\lf/M\cong\g$$
is a simple Lie algebra, we conclude that $\lf$ is Lie-simple.
\end{proof}

A left Leibniz algebra $\lf$ is called {\em simple\/} if $0$, $\leib(\lf)$,
$\lf$ are the only ideals of $\lf$, and $\leib(\lf)\subsetneqq\lf^2$,
(see \cite{AD}\footnote{Unfortunately, I have no access to \cite{AD}.
But in the paper \cite{DA}, which has been written by the same authors
around the same time, the non-triviality condition $\leib(\lf)\subsetneqq
\lf^2$ is omitted (see \cite[Definition 1]{DA}). (Note that the non-triviality
condition is also omitted in \cite[Definition 2.7]{AKOZ}.) But this
condition is needed to assure that the one-dimensional Lie algebra
and the two-dimensional solvable non-Lie Leibniz algebras are excluded,
and consequently, simple Leibniz algebras are semi-simple (see the
proof of \cite[Proposition 7.3]{F}).}, \cite[Definition 2.2]{ORT},
\cite[Section 1, p.\ 527]{GVKO}, and \cite[Definition 5.1]{DMS}).

Note that Theorems \ref{semisim} and \ref{liesim} are based on
the fact that a Leibniz algebra is semi-simple (resp.~Lie-simple) exactly
when its canonical Lie algebra is semi-simple (resp.~simple). The
analogous statement is not true for simple Leibniz algebras. But
somewhat surprisingly, in this case we can prove an even stronger
statement. Of course, it is not surprising that the proof of this result
is more involved than the proofs of the two previous results. It is
quite remarkable that this simplicity criterion neither needs a hypothesis
on the dimensions of the factors of the hemi-semidirect product nor
on the characteristic of the ground field.

\begin{thm}\label{sim}
Let $\g$ be a Lie algebra, and let $M$ be a left $\g$-module. Then
the left hemi-semidirect product $\g\ltimes_\ell M$ is simple if, and
only if, $\g$ is simple and either $M=0$ or $M$ is a non-trivial
irreducible $\g$-module.
\end{thm}

\begin{proof}
Let $\lf:=\g\ltimes_\ell M$ denote the left hemi-semidirect product
of $\g$ and $M$, and suppose that $\lf$ is simple. Firstly, we show
that $\leib(\lf)=M$. It follows from Lemma \ref{prophemisemi}\,(a)
that $\leib(\lf)=\g M\subseteq M$. Since $M$ is an ideal of $\lf$
and by hypothesis $\lf$ is simple, we have that $M=\leib(\lf)$ or
$M=\lf$. But in the second case we obtain that $\lf=M$ is an
abelian Lie algebra which contradicts the simplicity of $\lf$.

Moreover, from $\leib(\lf)=M$ and \cite[Proposition 7.2]{F}
we conclude that
$$\g\cong\lf/M=\lf/\leib(\lf)=\lf_\lie$$
is simple.

Secondly, we prove that $M$ is irreducible. So let $N$ be a
$\g$-submodule of $M=\leib(\lf)$. Since $N$ is an ideal of
$\lf$ and by hypothesis $\lf$ is simple, we have that $N=0$
or $N=\leib(\lf)=M$.

Thirdly, we show that either $M=0$ or $\g M\ne 0$. Suppose
to the contrary that $M\ne 0$ and $\g M=0$. But this implies
that $\lf=\g\times M$ is a Lie algebra, and thus either $M$ is
a non-zero proper ideal of $\lf$ or $\lf=M$ is abelian which
again is a contradiction.

Conversely, suppose now that $\g$ is simple and either $M=0$
or $M$ is a non-trivial irreducible $\g$-module. In the first case
we have that $\lf=\g$ is simple. So it remains to consider the
second case. Since $M$ is $\g$-irreducible and non-trivial,
we have that $\g M=M$, and thus it follows from Lemma
\ref{prophemisemi}\,(a) that $\leib(\lf)=\g M=M$.

Now consider the epimorphism $\pi:\lf\to\g$ of left Leibniz algebras
that is defined by $\pi(g,m):=g$ for any $g\in\g$ and $m\in M$.
Suppose that $\I$ is an ideal of $\lf$. Then $\pi(\I)$ is an ideal of
the simple Lie algebra $\g$, and therefore we obtain that $\pi(\I)
=0$ or $\pi(\I)=\g$.

In the case $\pi(\I)=0$ we have that $\I\subseteq M$. Since $M$
is an irreducible $\g$-module and $\I$ is an $\g$-submodule of $M$,
we conclude that $\I=0$ or $\I=M=\leib(\lf)$ as desired.

In the case $\pi(\I)=\g$ we have that $\lf=\I+M$. Namely, let
$(g,m)\in\lf=\g\ltimes_\ell M$ be arbitrary. Then it follows from
$\pi(\I)=\g$ that there exists an element $y\in\I$ such that
$\pi(y)=g$. Hence we obtain that $y=(g,m^\prime)$ for some
element $m^\prime\in M$, and therefore we conclude that
$$(g,m)=(g,m^\prime)+(0,m-m^\prime)=y+(0,m-m^\prime)
\in\I+M\,,$$
i.e., $\lf\subseteq\I+M$.

Since by hypothesis $M$ is an irreducible $\g$-module and
$\I\cap M$ is an $\g$-submodule of $M$, we have that either
$\I\cap M=0$ or $\I\cap M=M$. In the former case we deduce
that
$$\lf M=(\I+M)M=\I M\subseteq\I\cap M=0\,,$$
which cannot be true, as by thypothesis, the action of $\lf$
on $M$ is non-trivial. Consequently, we have that $\I\cap M=
M$, i.e., $M\subseteq\I$, and therefore the last inclusion implies
that $\lf=\I+M\subseteq\I$, i.e., $\I= \lf$.

Finally, we have to show that $\leib(\lf)\subsetneqq\lf^2$.
Suppose to the contrary that $\leib(\lf)=\lf^2$. Then we have
that $\g\cong\lf/M=\lf/\leib(\lf)=\lf/\lf^2$ is abelian which
contradicts the hypothesis that $\g$ is simple.
\end{proof}

\noindent {\bf Remark 1.} In order to show $\leib(\lf)=M$ in the
``only if"-part of the proof of Theorem~\ref{sim} one can also
argue as follows. Namely, we obtain from \cite[Propositions~7.3
and 7.5]{F} in conjunction with Lemma \ref{prophemisemi}\,(a)
and (c) that
$$M\subseteq C^\ell(\lf) =\leib(\lf)=\g M\subseteq M\,.$$
\vspace{-.3cm}

It might be possible that the ``if"-part of Theorem \ref{sim} has
already been obtained in \cite{AD}.\footnote{As mentioned earlier,
I have no access to \cite{AD}. What is stated here is taken from
the papers \cite{ORT} and \cite{GVKO}. On the other hand, in
the paper \cite{DA}, which has been written by Abdukassymova
and Dzhumadil’daev around the same time as \cite{AD}, the
``if"-part of Theorem \ref{sim} is not explicitly mentioned.
Moreover, on p.\ 1049 of \cite{DA} a related discussion does
not exclude the trivial irreducible module in the definition of a
{\em standard simple Leibniz algebra\/}, The latter should be
exactly the left hemi-semidirect products considered in Theorem
\ref{sim}.} This result is stated, although without the hypothesis
that $M$ is non-trivial and without any indication of a proof, in
Example 2.3 of \cite{ORT}\footnote{In \cite[Example 5.3]{DMS}
the authors refer to this result (and not to \cite{AD}) to justify
that the hemi-semidirect products $\slf_2(\C)\ltimes_\ell L$ are
simple for any irreducible left $\slf_2(\C)$-module $L$ of dimension
$\ge 2$.
Of course, the simplicity of these algebras is an immediate
consequence of the ``if"-part of Theorem \ref{sim}.} and
again in Example~1 at the end of Section 2 of \cite{GVKO}.
But clearly, the ``if"-part of Theorem \ref{sim} is not true
without assuming that $M$ is non-trivial.

Let $\C[t,t^{-1}]$ denote the algebra of Laurent polynomials
in one variable, let
$$\mathcal{W}:=\der(\C[t,t^{-1}])=\C[t,t^{-1}]\frac{d}{dt}
=\bigoplus_{n\in\Z}\C d_n$$
be the two-sided Witt algebra, where $d_n:=t^{n+1}\frac{d}{dt}$.
It is well known that $\mathcal{W}$ is a simple Lie algebra (for example,
see \cite[Theorem 1.2]{Pa}). Now let
$$V_{\alpha,\beta}=\bigoplus_{m\in\Z}\C v_m$$
be the left $\mathcal{W}$-module with action
$$d_n\cdot v_m:=[\alpha+\beta(n+1)+m]v_{m+n}$$
for any integers $m,n\in\Z$ and any complex numbers $\alpha,\beta$.
Then one can apply Theorem \ref{sim} in conjunction with \cite[Proposition
1.1]{KRR} to determine when the hemi-semidirect products $\mathcal{W}
\ltimes_\ell V_{\alpha,\beta}$ are simple (see \cite[Theorems 3.1 and
3.2]{COK} for the ``only if"-part).\footnote{In Remark 1 on p.\ 2062
of \cite{COK} the authors claim (without proof) that the Leibniz algebras
they consider are simple. On the other hand, more generally, they also
investigate non-split Leibniz extensions of $\mathcal{W}$ by the anti-symmetrization
of $V_{\alpha,\beta}$ which we will only discuss in the second part of this paper
(see also \cite{D1} and \cite{DA} for non-split Leibniz extensions of
simple Lie algebras by irreducible anti-symmetric Leibniz bimodules).}

\begin{pro}
The hemi-semi-direct product $\mathcal{W}\ltimes_\ell V_{\alpha,
\beta}$ is simple if, and only if, $\alpha\not\in\Z$ or $\alpha\in
\Z$ and $\beta\ne 0,1$.
\end{pro}



\section{Structure theorems}\label{structhm}


As a consequence of Theorem \ref{semisim} we obtain the
following characterization of finite-dimensional semi-simple left
Leibniz algebras over fields of characteristic zero. Note that
the implication (i)$\Rightarrow$(ii) in Theorem \ref{semisimchar0}
is \cite[Theorem~3.4\,(a) \& (b)]{AKOZ} as the hypothesis
that $\lf$ is indecomposable is not used in this part of the
theorem (see also \cite[Proposition 2.11]{AKOZ}). Moreover,
note that Theorem \ref{semisimchar0} is a generalization of
the first and second statement of \cite[Theorem 7.12]{F}
and compare part of the argument also with the proof of
\cite[Lemma 4.1]{FW}).

\begin{thm}\label{semisimchar0}
Let $\lf$ be a left Leibniz algebra over a field $\F$ of characteristic
zero. Then the following statements are equivalent:
\begin{enumerate}
\item[{\rm (i)}] $\lf$ is finite-dimensional semi-simple.
\item[{\rm (ii)}] There exists a finite-dimensional  semi-simple Lie
                          subalgebra $\g$ of $\lf$ such that $\lf=\g\oplus
                          \leib(\lf)$, and $\leib(\lf)$ is a finite-dimensional
                          completely reducible anti-symmetric $\g$-bimodule
                          that does not contain any non-zero trivial
                          $\g$-subbimodule.
\item[{\rm (iii)}] $\lf$ is isomorphic to the left hemi-semidirect
                          product of a finite-dimensional semi-simple Lie
                          algebra $\g$ over $\F$ and a finite-dimensional
                          completely reducible left $\g$-module $M$ that
                          does not contain any non-zero trivial $\g$-sub\-module.
\end{enumerate}
\end{thm}

\begin{proof}
We begin by proving the implication (i)$\Rightarrow$(ii). According
to Levi's theorem for Leibniz algebras (see \cite[Proposition 2.4]{P}
or \cite[Theorem 1]{B}), there exists a semi-simple Lie subalgebra
$\g$ of $\lf$ such that $\lf=\g\oplus\rad(\lf)$. Since by hypothesis
$\lf$ is semi-simple, we have that $\rad(\lf)=\leib(\lf)$, and thus
we obtain that $\lf=\g\oplus\leib(\lf)$. Moreover, we conclude
from the example in Section~\ref{prelim} that $M:=\leib(\lf)$
is an anti-symmetric $\g$-bimodule. We have either $M=0$ or
$M\ne 0$. In the former case the assertion is clear, and so we
assume in the remainder of the proof of this implication that
$M\ne 0$. Because $M$ is an anti-symmetric $\g$-bimodule,
being an $\g$-subbimodule is the same as being a left $\g$-submodule.
Hence Weyl's theorem on complete reducibility (see \cite[Theorem
6.3]{H}) implies that $M$ is a completely reducible left $\g$-module
as well as a completely reducible $\g$-bimodule. Finally, let $M_\triv$
denote the largest trivial $\g$-sub(bi)module of $M$. As $M$ is
completely reducible, there exists an $\g$-sub(bi)module $C$ of
$M$ such that $M=M_\triv\oplus C$. Since $\g$ is a Lie algebra
and $M\lf=0$, we obtain that $(g+m)(g+m)=g\cdot m$ for any
elements $g\in\g$ and $m\in M$. This shows that $M=\leib(\lf)=
\g M$, and therefore we obtain
$$M_\triv\subseteq M=\g M=\g M_\triv+\g C=\g C\subseteq C\,,$$
i.e., $M_\triv=M_\triv\cap C=0$. Hence $M$ does no contain any
non-zero trivial $\g$-sub(bi)module.

Next, we prove the implication (ii)$\Rightarrow$(iii). Set $M:=\leib
(\lf)$. Again in the case of $M=0$ there is nothing to prove. So we
assume that $M\ne 0$. Then we have to show that $\lf=\g\oplus M$
is a left hemi-semidirect product. Since $M=\leib(\lf)$ is an abelian
ideal of $\lf$ such that $M\lf=0$, we obtain that
$$(g_1+m_1)(g_2+m_2)=[g_1,g_2]+g_1\cdot m_2+m_1\cdot g_2=
[g_1,g_2]+g_1\cdot m_2$$
for any $g_1,g_2\in\g$ and $m_1,m_2\in M$, and therefore the
function $g+m\mapsto(g,m)$ defines an isomorphism of left
Leibniz algebras from $\lf$ onto the left hemi-semidirect product
$\g\ltimes_\ell M$ of $\g$ and $M$.

Finally, we prove the remaining implication (iii)$\Rightarrow$(i).
By hypothesis, we have that $\g S=S$ for every irreducible direct
summand $S$ of M, and therefore the same holds for $M$, namely,
$\g M=M$, and thus the assertion is an immediate consequence
of Theorem \ref{semisim}.
\end{proof}

\noindent {\bf Remark 2.} More precisely than part (iii) of Theorem
\ref{semisimchar0}, Theorem 3.4\,(c) in \cite{AKOZ}
gives a description of the irreducible summands of $M$ as modules
over the simple factors of $\g$.
\vspace{.2cm}

\noindent {\bf Remark 3.} Let $\rk(\lf):=\dim_\F\leib(\lf)$ denote
the {\em rank\/} of a Leibniz algebra $\lf$ over a field $\F$ (see
\cite[Section 5, p.\ 8]{MY}). Then it follows from Theorem
\ref{semisimchar0} that for every finite-dimensional semisimple
Leibniz algebra $\lf$ over a field of characteristic zero we have
that $\rk(\lf)\ne 1$.
\vspace{.2cm}

Next, we characterize Lie-simple Leibniz algebras over fields
of characteristic zero in a similar way as we did for semi-simple
Leibniz algebras in the previous result. The proof of Theorem
\ref{liesimchar0} is very similar to the proof of Theorem
\ref{semisimchar0} and is left to the reader.

\begin{thm}\label{liesimchar0}
Let $\lf$ be a left Leibniz algebra over a field $\F$ of characteristic
zero. Then the following statements are equivalent:
\begin{enumerate}
\item[{\rm (i)}] $\lf$ is finite-dimensional Lie-simple.
\item[{\rm (ii)}] There exists a finite-dimensional simple Lie subalgebra
                         $\g$ of $\lf$ such that $\lf=\g\oplus\leib(\lf)$, and
                         $\leib(\lf)$ is a finite-dimensional completely reducible
                         anti-symmetric $\g$-bimodule that does not contain
                         any non-zero trivial $\g$-subbimodule.
\item[{\rm (iii)}] $\lf$ is isomorphic to the left hemi-semidirect product
                          of a finite-dimensional simple Lie algebra $\g$ over $\F$
                          and a finite-dimensional completely reducible left $\g$-module
                          $M$ that does not contain any non-zero trivial $\g$-submodule.
\end{enumerate}
\end{thm}

Note that \cite[Theorem 3.5]{ORT} is an immediate consequence of
Theorem \ref{liesimchar0}, Weyl's theorem on complete reducibility,
and the well-known structure of finite-dimensional irreducible $\slf_2
(\C)$-modules (see \cite[Theorem 7.2]{H}).

Finally, we characterize finite-dimensional  simple Leibniz algebras
over fields of characteristic zero. Note that Theorem \ref{simchar0}
is a generalization of the third statement of \cite[Theorem 7.12]{F}.

\begin{thm}\label{simchar0}
Let $\lf$ be a left Leibniz algebra over a field $\F$ of characteristic
zero. Then the following statements are equivalent:
\begin{enumerate}
\item[{\rm (i)}] $\lf$ is finite-dimensional simple.
\item[{\rm (ii)}] There exists a finite-dimensional simple Lie subalgebra
                         $\g$ of $\lf$ such that $\lf=\g\oplus\leib(\lf)$, and either
                         $\leib(\lf)=0$ or $\leib(\lf)$ is a non-trivial finite-dimensional
                          irreducible anti-symmetric $\g$-bimodule.
\item[{\rm (iii)}] $\lf$ is a finite-dimensional simple Lie algebra or $\lf$
                          is isomorphic to the left hemi-semidirect product of a
                          finite-dimensional simple Lie algebra $\g$ over $\F$
                          and a non-trivial finite-dimensional irreducible left
                          $\g$-module $M$.
\end{enumerate}
\end{thm}

\begin{proof}
Since the implication (ii)$\Rightarrow$(iii) follows as in the proof of
Theorem \ref{semisimchar0} and the implication (iii)$\Rightarrow$(i)
is an immediate consequence of Theorem \ref{sim}, it remains to
prove the implication (i)$\Rightarrow$(ii). Since every simple Leibniz
algebra is semi-simple (see \cite[Proposition 7.3]{F}), the assertion
can be obtained from Theorem \ref{semisimchar0} and the third
statement of \cite[Theorem 7.12]{F}.
\end{proof}

In the first paragraph of p.\ 2009 in \cite{AKOZ} the authors
claim (without any further justification) that the implication 
(i)$\Rightarrow$(iii) of Theorem \ref{sim} follows from
\cite[Proposition~2.11]{AKOZ}\footnote{It is remarkable
that the non-triviality condition on the Leibniz kernel appears in
\cite[Proposition~2.11]{AKOZ}, but not in the statement on
the structure of simple Leibniz algebras which is supposedly
deduced from it.}

Note that \cite[Theorem 3.4]{ORT} is an immediate consequence
of Theorem \ref{simchar0} and the well-known structure of
finite-dimensional irreducible $\slf_2(\C)$-modules. In addition,
Theorem \ref{simchar0} shows the converse, namely, the Leibniz
algebras in \cite[Theorem~3.4]{ORT} are simple if, and only if, the
Leibniz kernel $I$ is not isomorphic to the trivial irreducible $\slf_2
(\C)$-(bi)module.\footnote{It seems this is what the authors
intended to prove (see Example 2.3 and the last paragraph on
p.\ 1510 of \cite{ORT}). As one can see from Theorem 3.1 in
the later paper \cite{GVKO}, the authors still only address the
``only if"- part (although omitting the non-triviality condition
on the Leibniz kernel).}

It is well known (see \cite[p.\ 529]{GVKO} and \cite[Propositions 7.2
and 7.8]{F}) that the following inclusions hold:
$$\{\,\mbox{simple}\,\}\subsetneqq\{\,\mbox{Lie-simple}\,\}
\subsetneqq\{\,\mbox{semi-simple}\,\}\,.$$
Moreover, it is immediately clear from the structure theorems
in this section that every inclusion is strict.

A classical result says that every finite-dimensional Lie algebra
over a field of characteristic zero is semi-simple exactly when
it is the direct product of simple Lie algebras (see \cite[Theorem
5.2]{H} for the more substantial implication). One implication is
still true for Leibniz algebras, namely, that a finite direct product
of semi-simple Leibniz algebras is semi-simple, but the converse
does not hold as has been shown by Example 2 in \cite{GVKO}
(see also \cite[Examples 2.12 and 3.6]{AKOZ}). More generally,
in Example 3 of \cite{GVKO} it has been demonstrated
that a finite-dimensional semi-simple Leibniz algebra is not necessarily
a direct product of Lie-simple Leibniz algebras.\footnote{At the end
of \cite{GVKO} the authors (without any proof) give a characterization
of those semi-simple Leibniz algebras that are a direct sum of Lie-simple
Leibniz algebras.} By using Theorems \ref{semisimchar0} and Theorem
\ref{liesimchar0} we give a general construction to establish the last
statement. For this we will need the following description of the Leibniz
kernel of a finite direct product of Leibniz algebras.

\begin{lem}\label{leibnizker}
Let $\lf_1$, ..., $\lf_n$ be Leibniz algebras. Then
$$\leib(\lf_1\times\cdots\times\lf_n)=\leib(\lf_1)\times\cdots
\times\leib(\lf_n)\,.$$
\end{lem}

\begin{proof}
Let $x=(x_1,\dots,x_n)\in\lf:=\lf_1\times\cdots\times\lf_n$ be
arbitrary. Then
$$x^2=(x_1^2,\dots,x_n^2)\in\leib(\lf_1)\times\cdots\times
\leib(\lf_n)\,,$$
which implies that $\leib(\lf):=\langle x^2\mid x\in\lf\rangle_\F
\subseteq\leib(\lf_1)\times\cdots\times\leib(\lf_n)$.

Next, we prove the reverse inclusion. For this we let the element
$x=(x_1,\dots,x_n)\in\leib(\lf_1)\times\cdots\times\leib(\lf_n)$ be
arbitrary. Then for every $j\in\{1,\dots,n\}$ there exist $\alpha_{ji_j}
\in\F$ and $x_{ji_j}\in\lf_j$ such that
$$x_j=\sum_{i_j= 1}^{m_j}\alpha_{ji_j}x_{ji_j}^2\,,$$
and it follows that
\begin{eqnarray*}
x & = & (x_j)_{1\le j\le n}=\left(\sum_{i_j=1}^{m_j}\alpha_{ji_j}
x_{ji_j}^2\right)_{1\le j\le n}=\sum_{j=1}^n\sum_{i_j=1}^{m_j}
\alpha_{ji_j}\left(x_{ji_j}^2\right)_{1\le j\le n}\\
& = & \sum_{j=1}^n\sum_{i_j=1}^{m_j}\alpha_{ji_j}\left(x_{ji_j}
\right)_{1\le j\le n}^2\in\leib(\lf)\,.
\end{eqnarray*}
This completes the proof.
\end{proof}

\noindent {\bf Remark 4.} Since we did not use the Leibniz identity
in the proof of Lemma \ref{leibnizker} a similar statement holds for
the subspace spanned by the squares in a finite product of arbitrary
algebras.
\vspace{.2cm}

Now we are ready to give the following general construction of
semi-simple left Leibniz algebras that are not a direct product of
Lie-simple Leibniz algebras. Note that Proposition \ref{countexam}
conceptualizes and simplifies Example 3 in \cite{GVKO} and
Example~2.12 in \cite{AKOZ} considerably.

\begin{pro}\label{countexam}
Let $\g_1$ and $\g_2$ be two (not necessarily distinct)
finite-dimensional simple Lie algebras over a field of characteristic
zero, and let $M$ be a finite-dimensio\-nal non-trivial irreducible
left $(\g_1\times\g_2)$-module. Then $\lf:=(\g_1\times\g_2)
\ltimes_\ell M$ is a semi-simple left Leibniz algebra that is not
isomorphic to a direct product of Lie-simple Leibniz algebras.
\end{pro}

\begin{proof}
According to Theorem \ref{semisimchar0}, we have that $\lf$ is
semi-simple. Suppose that $\lf=\lf_1\times\cdots\times\lf_n$
is a finite direct product of Lie-simple Leibniz algebras $\lf_1$,
\dots, $\lf_n$. (Since $\lf$ is finite-dimensional, the direct
product necessarily is finite.) Then it follows from the proof
of Theorem \ref{semisim} in conjunction with Lemma
\ref{leibnizker} that $M=\leib(\lf_1)\times\cdots\times\leib
(\lf_n)$, and therefore we obtain from the irreducibility of
$M$ that $n=1$. Hence, $\lf=\lf_1$ is Lie-simple. But it is a
consequence of Theorem~\ref{liesimchar0} that $\lf$ is not
Lie-simple as its canonical Lie algebra $\g_1\times\g_2$ is not
simple, which is a contradiction. 
\end{proof}

\noindent {\bf Remark 5.} The construction in Proposition
\ref{countexam} is in a sense dual to the one in \cite[Example
7.13]{F} in which a Lie-simple Leibniz algebra was given that is
not simple (see also the second part of \cite[Example 5.3]{DMS}).
\vspace{.2cm}


\section{Applications}\label{appl}


In this last section we give some applications of our structure
theorem for finite-dimensional semi-simple Leibniz algebras in
characteristic zero (see Theorem \ref{semisimchar0}).

In fact, the next result just needs Levi's theorem for Leibniz
algebras (see the proof of the implication (i)$\Rightarrow$(ii)
of Theorem \ref{semisimchar0}).

\begin{pro}\label{grad}
Every finite-dimensional semi-simple left Leibniz algebra over
a field of characteristic zero is $\Z_2$-graded.
\end{pro}

\begin{proof}
Let $\lf=\s\oplus\leib(\lf)$ be a finite-dimensional semi-simple
left Leibniz algebra over a field of characteristic zero with
Levi subalgebra $\s$. Set $\lf_{\overline{0}}:=\s$ and
$\lf_{\overline{1}}:=\leib(\lf)$. Since $\s$ is a subalgebra,
we have that $\lf_{\overline{0}}\lf_{\overline{0}}=\s\s
\subseteq\s=\lf_{\overline{0}}$. Moreover, as $\leib(\lf)$
is an ideal of $\lf$, we obtain that $\lf_{\overline{0}}\lf_{
\overline{1}}=\s\leib(\lf)\subseteq\leib(\lf)=\lf_{\overline{1}}$
and $\lf_{\overline{1}}\lf_{\overline{0}}=\leib(\lf)\s\subseteq
\leib(\lf)=\lf_{\overline{1}}$. (In fact, we have that $\leib(\lf)
\s=0$.) Finally, as $\leib(\lf)$ is abelian, we have that $\lf_{
\overline{1}} \lf_{\overline{1}}=0\subseteq\lf_{\overline{0}}$.
\end{proof}

\noindent {\bf Remark 6.} In fact, slightly more general, the
proof of Proposition \ref{grad} shows that every finite-dimensional
left Leibniz algebra with an abelian radical over a field of
characteristic zero is $\Z_2$-graded.
\vspace{.2cm}

Furthermore, the proof of Proposition \ref{grad} shows that
every finite-dimensional semi-simple left Leibniz algebra is a left
Leibniz superalgebra. Here a {\em left Leibniz superalgebra\/}
is a superalgebra for which every left multiplication operator is a
superderivation (see \cite[Definitions 1.1.4 and 1.1.2]{K} and
\cite{D2}).

Theorem \ref{semisimchar0} in conjunction with Lemma
\ref{prophemisemi}\,(d) implies the vanishing of the right
center of a finite-dimensional semi-simple left Leibniz
algebra in characteristic zero.

\begin{pro}\label{rightcent}
If $\lf$ is a finite-dimensional semi-simple left Leibniz algebra
over a field of characteristic zero, then $C^r(\lf)=0$.
\end{pro}

\begin{proof}
It follows from Theorem \ref{semisimchar0} that $\lf$ is the
left hemi-semidirect product of a finite-dimensional semi-simple
Lie algebra $\g$ and a finite-dimensional completely reducible
left $\g$-module $M$ that does not contain any non-zero trivial
$\g$-submodule. In particular, we have that $C(\g)=0$ and
$M^\g=0$. Now apply Lemma \ref{prophemisemi}\,(d).
\end{proof}

Recently, Friedrich Wagemann observed that the two-dimensional
nilpotent non-Lie Leibniz algebra $\nf=\F e\oplus\F c$
with $ee=ec=ce=0$ and $cc=e$ is a Lie superalgebra if $e$ has
even degree and $c$ has odd degree. In fact, the Lie superalgebra
$\nf$ has already been considered by Kac, namely, it is the
Heisenberg superalgebra $N^\prime$ for $n=0$ (see \cite[Example
on p.\ 18 in Section 1.1.6]{K}). On the other hand, although
being $\Z_2$-graded by Proposition \ref{grad}, a finite-dimensional
semi-simple left non-Lie Leibniz algebra is never a Lie superalgebra.
(Otherwise the super-anticommutativity would imply that the left
and right center are the same which according to \cite[Proposition
7.5]{F} and Proposition \ref{rightcent} cannot happen.)

Following Mason and Yamskulna \cite[p.\ 2]{MY} we call a
left Leibniz algebra $\lf$ {\em left central\/} if $\leib(\lf)
\subseteq C^r(\lf)$. They observed that the following
inclusions are strict:
$$\{\,\mbox{Lie}\,\}\subsetneqq\{\,\mbox{symmetric Leibniz}
\,\}\subsetneqq\{\,\mbox{left central Leibniz}\,\}\subsetneqq
\{\,\mbox{left Leibniz}\,\}\,.$$


For example, note that it follows from Lemma \ref{prophemisemi}\,(a)
and (d) that the two-dimensional solvable left Leibniz algebra $\af=\F
h\ltimes_\ell\F e$ is not left central (see also \cite[Examples~2.10
and 2.16]{F}) which shows that the third inclusion is strict.

On the other hand, it is a consequence of Proposition \ref{rightcent}
that, in the above chain of inclusions considered only for finite-dimensional
semi-simple Leibniz algebras over fields of characteristic zero, the
first two inclusions are equalities:

\begin{cor}\label{symmsemisim}
Every finite-dimensional semi-simple left central Leibniz algebra
over a field of characteristic zero is a Lie algebra.
\end{cor}

Since there is a plethora of finite-dimensional semi-simple left non-Lie
Leibniz algebras over a field of characteristic zero (see Theorem
\ref{semisimchar0}), Corollary \ref{symmsemisim} yields a multitude
of left Leibniz algebras that are not left central.
\vspace{.2cm}

\noindent {\bf Remark 7.} In fact, much more is true. Namely,
one can deduce from Lemma~\ref{prophemisemi}\,(a) and (d)
the following statement:
If $\g$ is perfect and $M$ is a non-trivial $\g$-module, then the
hemi-semidirect product $\g\ltimes_\ell M$ is not left central.\footnote{Kinyon and Weinstein
already observed that $\g\ltimes_\ell M$ is not a Lie algebra when $\g M\ne 0$
(see \cite[p.\ 530]{KW}).}
\vspace{.2cm}

We finish this section by applying Theorem \ref{semisimchar0}
in conjunction with \cite[Theorem~4.5\,(c)]{FW} to generalize
\cite[Theorem 4.5\,(c)]{AKOZ} from the complex numbers to arbitrary
fields of characteristic zero.

\begin{thm}\label{semisimder}
Let $\lf=\s\oplus\leib(\lf)$ be a finite-dimensional semi-simple
Leibniz algebra over a field $\F$ of characteristic zero with Levi
factor $\s$. Then $$\der(\lf)\cong\s\oplus\bigoplus_{i=1}^s
\bigoplus_{j=1}^t\Hom_\s(\s_i,M_j)^{\oplus a_ib_j} \oplus
\bigoplus_{k=1}^t\End_\s(M_k)^{\oplus b_k^2}$$ as $\F$-vector
spaces, where $\s=\bigoplus\limits_{i=1}^s\s_i^{\oplus a_i}$
is the decomposition of $\s$ into a direct sum of pairwise
non-isomorphic simple ideals and $\leib(\lf)=\bigoplus\limits_{j=1}^t
M_j^{\oplus b_j}$ is the decomposition of $\leib(\lf)$ into a
direct sum of pairwise non-isomorphic non-trivial irreducible
left $\s$-submodules.
\end{thm}

\begin{proof}
According to Theorem \ref{semisimchar0}, we have that $\lf=
\s\oplus\leib(\lf)$, where $\s$ is a semi-simple Lie subalgebra
of $\lf$ and $\leib(\lf)$ is a completely reducible anti-symmetric
$\s$-bimodule that does not contain any non-zero trivial
$\s$-subbimodule. As $M:=\leib(\lf)$ is anti-symmetric, the
latter is the same as saying that $\leib(\lf)$ is a completely
reducible left $\s$-module that does not contain any non-zero
trivial $\s$-submodule. So we obtain that $M=\bigoplus\limits_{j=1}^t
M_j^{\oplus b_j}$ is a direct sum of pairwise non-isomorphic
non-trivial irreducible $\s$-submodules. Moreover, it follows
from \cite[Theorem 5.2]{H} that $\s$ is a direct sum of ideals
that are simple as Lie algebras. Observe that every such ideal
of $\s$ is an irreducible left $\s$-module. Hence we have that
$\s=\bigoplus\limits_{i=1}^s\s_i^{\oplus a_i}$ is a direct sum
of pairwise non-isomorphic irreducible left $\s$-modules.

Note that $L(\lf):=\{L_x\mid x\in\lf\}$ is the set of inner derivations
of $\lf$ (see \cite[Section~4]{F}), and therefore we obtain from
\cite[Proposition 4.3]{F} and \cite[Theorem~4.5\,(c)]{FW} that
as $\F$-vector spaces $$\der(\lf)\cong L(\lf)\oplus\Hom_\lf
(\lf_{\ad,\ell},M)\,,$$ where $\lf_{\ad,\ell}$ denotes the left
adjoint $\lf$-module. Furthermore, we deduce from
\cite[Proposition 2.23]{F} and \cite[Proposition 7.5]{F} that
$$L(\lf)\cong\lf/C^\ell(\lf)=\lf/\leib(\lf)\cong\s$$ as Lie algebras,
and thus as $\F$-vector spaces.

Let us now consider the second summand $\Hom_\lf(\lf_{\ad,\ell},
M)$. By virtue of \cite[Lemma 3.3]{F}, we have $\leib(\lf)M=0$,
and thus we conclude that
$$\Hom_\lf(\lf_{\ad,\ell},M)=\Hom_\s(\lf_{\ad,\ell},M)\,.$$
As left $\s$-modules, we have that
$$\lf_{\ad,\ell}\cong\s\oplus M=\bigoplus_{i=1}^s
\s_i^{\oplus a_i}\oplus\bigoplus_{j=1}^tM_j^{\oplus b_j}\,.$$
Hence, we obtain from Schur's lemma that
$$\Hom_\s(\lf_{\ad,\ell},M)\cong\bigoplus_{i=1}^s
\bigoplus_{j=1}^t\Hom_\s(\s_i,M_j)^{\oplus a_ib_j}\oplus
\bigoplus_{k=1}^t\End_\s(M_k)^{\oplus b_k^2}$$
as $\F$-vector spaces, which completes the proof.
\end{proof}

As an immediate consequence of (the proof of) Theorem
\ref{semisimder} we obtain a new very short proof of
\cite[Corollary 4.6]{FW}.

\begin{cor}
Every finite-dimensional semi-simple non-Lie Leibniz algebra
over a field of characteristic zero has derivations that are not inner.
\end{cor}

\begin{proof}
If $\lf$ is a non-Lie Leibniz algebra, then using the notation of
Theorem \ref{semisimder} we have that $t\ge 1$ and $M_1\ne 0$.
Therefore, $0\ne\id_{M_1}\in\Hom_\s(\lf_{\ad,\ell},M)$ implies
that $\der(\lf)/L(\lf)\cong\Hom_\s(\lf_{\ad,\ell},M)\ne 0$.
\end{proof}


\noindent {\bf Acknowledgment.} 
The author would like to thank Friedrich Wagemann for sharing his
observations on the two-dimensional non-Lie algebra $\nf$ of Heisenberg
type which are mentioned in the paragraph after the proof of Proposition
\ref{rightcent}.



\end{document}